\documentclass{article} 

\usepackage{graphicx}

\usepackage{amsmath}
\usepackage{amsfonts}
\usepackage{subfig}
\usepackage{float}
\usepackage{hyperref}

\usepackage{algorithm}
\usepackage{algpseudocode}

\makeatletter
\algrenewcommand\ALG@beginalgorithmic{\footnotesize}
\makeatother

\newcommand{\F}{\mathcal{F}}
\newcommand{\T}{\mathcal{T}}

\newcommand{\QQ}{\mathbf{Q}}
\newcommand{\ZZ}{\mathbf{Z}}

\newcommand{\CC}{\mathbb C}
\newcommand{\HH}{{\mathcal{H}}}
\newcommand{\DD}{{\mathbb{D}}}
\newcommand{\FT}{\mathcal{FT}}

\newcommand{\psl}{\mathrm{PSL}_2 (\ZZ)}

\newcommand{\ol}[1]{\overline{#1}}
\newcommand{\cark}{\c{c}ark}
\newcommand{\tr}{\mathrm{Tr}}
\newcommand{\im}{\mathrm{Im}}
\newcommand{\oo}{\mathcal{O}}
\newcommand{\lra}{\longrightarrow}
\renewcommand{\SS}{\begin{pmatrix}
                 0 & -1\\ 1 & 0
                 \end{pmatrix}
}
\newcommand{\LL}{\begin{pmatrix}
                 1 & -1\\ 1 & 0
                 \end{pmatrix}
}
\newcommand{\dd}{\mathcal{D}}

\title{InfoMod: A visual and computational approach to Gauss' binary quadratic forms}
\author{Ayberk Zeytin\footnotemark[1], Hakan Ayral\footnotemark[1], A. Muhammed Uluda\u{g}\footnotemark[1]}

\begin{document}

\maketitle

\begin{abstract}
InfoMod is a new software and application devoted to the modular group, $\psl$. It has algorithms that deals with the classical correspondences among continued fractions, geodesics on the modular surface and binary quadratic forms. In addition the software implements the recently discovered representation of Gauss' indefinite binary quadratic forms and their classes in terms of certain infinite planar graphs (dessins) called \c carks. InfoMod illustrates various aspects of these forms, i.e. Gauss' reduction algorithm, the representation problem of forms, ambiguous and reciprocal forms. It can be used as an educational tool, and might be used to explore some new facts about these objects. 
\end{abstract}

\footnotetext[1]{
Galatasaray University, \c Ciragan Cad. No 36 Be\c siktas\c s 34349  Istanbul, Turkey \\
This work is supported by the T\"UB\.ITAK grant 113R017 and the GSU grant 15.504.002.
}

\section{Introduction}

InfoMod is an innovative software which aims to visualize the deep arithmetic questions encircled by finitely generated infinite index subgroups of the modular group, which is by definition the the group consisting of $2\times2$ integral matrices of determinant $1$, denoted by $\psl$. These lead to the study of algebraic number theory of real quadratic number fields. Its intended use is instructive as well as research-level visual explorations of the modular group, binary quadratic forms, their geodesics, form classes, class groups and some algorithms such as Gauss' reduction algorithm.

One of the main problems of algebraic number theory is to understand factorization in the ring of integers of number field. The class group and its size, called the class number, are the prime indicators of how unique factorization in a ring of integer is. There are several important questions and conjectures about class groups and class numbers. One needs to recourse to algorithms for explicit computations of these groups. Since these problems resisted centuries of attacks, these are also important from the computational perspective and present many challenges. 

The simplest case of {\it quadratic} number fields, can be understood in terms of the binary quadratic forms of Gauss. Here the main difficulty is in the {\it real} quadratic case, which corresponds to indefinite binary quadratic forms. These forms can be represented as the fixed points of the elements of the modular group, and this group is the main hero of our paper.

The modular group is isomorphic to the free product $\ZZ/2\ZZ \ast \ZZ/3\ZZ$. As such its elements can be represented as the set of half-edges of an infinite rooted bipartite tree $\F\T$ with a planar structure, called the {\it Farey tree}. In order to visualize the elements of $\F\T$, we broaden the branches of $\F\T$ until they touch each other. The resulting mantra-like model of $\psl$ with cells representing its elements, fills the entire plane and is called the {\it sunburst}. It allows one to simultaneously represent up to approx. 1000 elements of the modular group on a computer screen with a standard resolution. By moving the centre of this sunburst (operation which is implemented in InfoMod), it is possible to represent another 1000-element portion of the modular group. By clicking over a cell, it is possible to obtain a detailed information on the group element represented by this cell. In particular, it is possible to see the corresponding binary quadratic form as well as its discriminant.

Now we alter our perspective slightly to get another representation of forms.
Suppose that the element $M\in\psl$ represents a reduced real binary quadratic form $f_M$. 
As a matter of fact, every element different from identity in the subgroup $\langle M\rangle$ also represents the same reduced form. 
The group $\langle M\rangle$ acts on the tree $\F$ in a natural way, and the quotient graph 
$\F\T/\langle M\rangle$ also represents $f_M$. Given an edge of the quotient graph, it is possible to recover the form.
These graphs are called {\it \c carks}. They look same as $\F\T$ except that they have exactly one circuit inside, called the {\it spine}. InfoMod makes it possible to pass, from the cell of $M$, to its \c cark. 

Many notions related to forms acquires a new formulation in terms of these graphs. For example, equivalent forms have isomorphic {\c carks}. 
The reduction algorithm of Gauss simply moves the edge characterizing the form towards the spine of the \c cark. This new formulation already served the first named author to give a solution of the representation problem of forms and improve on the Gauss reduction algorithm, \cite{reduction}. InfoMod makes it  possible to observe how these algorithms are carried out on the graph. 

Paper is organized as follows: In the next section we overview the theoretical background outlined in the above lines. More details can be found in \cite{UZD}. The third section is devoted to the software implementation of InfoMod, where we present a code library which help represent and do computations on binary quadratic forms, including reductions and enumeration of other forms which share a cycle with a reduced form. Later we present an interactive visualization application to explore the mentioned sunburst and \c carks visually. We plan to build into InfoMod some new features pertaining to the outer automorphism of the group PGL$(2,\mathbf Z)$ in the future. 

\section{Preliminaries}
\label{sec:preliminaries}
This section is devoted to recalling the basic facts around the modular group, binary quadratic forms, reduction theory and \c{c}arks. There are many beautiful texts on the subject among which we must note the first historically systematic treatment by Gauss, \cite{disquisitiones}. We refer to \cite{bqf/vollmer,buell/bqf} for a modern treatment of the topics related to binary quadratic forms.

\subsection{Modular group, binary quadratic forms and geodesics}
  
The modular group, $\psl$, acts on the upper half plane $\HH = \{z \in \CC \colon \im(z)>0\}$ via M\"{o}bius transformations, i.e. for $W = \begin{pmatrix}p & q \\r&s \end{pmatrix} \in \psl$ and $z \in \HH$, we have $W\cdot z := \frac{pz+q}{rz+s}$. This action leaves the Poincar\'{e} metric on $\HH$ invariant. Hence the geodesics of this metric, i.e. lines that are vertical to the real line and half circles whose centers are on the real line, are mapped onto geodesics by $\psl$. 

An \emph{integral binary quadratic form} (or BQF for short) is a homogeneous polynomial of degree two in two variables with integer coefficients. A BQF $f(x,y) = ax^{2}+bxy + cy^{2}$, is usually denoted by $(a,b,c)$ or its matrix form $\begin{pmatrix}
a & b/2 \\ b/2 & c                                                                                                                                                                                                                                                                                                                                                                                                                                                                                                                                                                                                                                                                                                                                              \end{pmatrix}$. We call $f$ primitive if the greatest common divisor of $a$, $b$ and $c$ is $1$. 

An element $W = \begin{pmatrix}p & q \\r&s \end{pmatrix} \in \psl$ has two fixed points, which are, by definition, roots of the quadratic equation $rz^{2} + (s-p)z -q = 0$. To $W$ we associate the BQF  $f_{W} := \frac{1}{\delta}\left(rx^{2} + (s-p)xy -qy^{2}\right)$; where $\delta$ is the greatest common divisor of $r$, $s-p$ and $q$, so that $f_{W}$ is by definition primitive. The discriminant of a BQF $f = (a,b,c)$ is defined as $\Delta(f) = b^{2}-4ac$. A BQF is called \emph{degenerate} if $\Delta(f)$ is a perfect square, \emph{positive definite} if $\Delta(f)<0$ and $a>0$, \emph{negative definite} if $\Delta(f)<0$ and $a<0$. Finally, a BQF is called \emph{indefinite} if $\Delta(f)>0$. 

We may classify elements on $\psl$ according to their traces. Namely, $W$ is called \emph{elliptic}, \emph{parabolic} and \emph{hyperbolic} if $|\tr(W)|<2$, $|\tr(W)| = 2$ and $|\tr(W)>2$, respectively. In other words, parabolic elements have a unique real fixed point. Elliptic elements have two fixed points one in the upper half plane and the other in the lower half plane. Hyperbolic elements have two distinct real fixed points. The geodesic in $\HH$ joining the two fixed points of a hyperbolic element $W \in \psl$\footnote{One may define attracting and repelling fixed points of $W$ by looking at its action along this geodesic. This distinction is unnecessary for our purposes.} which is merely a half circle of center ${(p-s)}/{2r}$ and radius $\frac{1}{2r}{\sqrt{\tr(W)^{2} - 4}}$ is called the geodesic associated to $W$ and denoted by $\gamma_{W}$. Observe that the classification of elements of $\psl$ and BQFs are quite parallel. Namely, parabolic elements give rise to degenerate BQFs, elliptic elements give rise to definite BQFs and hyperbolic elements give rise to indefinite BQFs.

\subsection{Reduction Theory of BQFs}

The modular group acts on the set $\mathcal{F}$ of non-degenerate binary quadratic forms by change of variables. More precisely, given $f = (a,b,c)$ and an element 
$M_{o} = 
	\begin{pmatrix} 
		p&q \\ r&s 
	\end{pmatrix}
$ we define $M_{o} \cdot f$ to be the BQF associated to the symmetric matrix $M_{o}^{t}M_{f}M_{o}$; where $M_{o}^{t}$ stands for the transpose of $M_{o}$. In particular, we say that two forms $f = (a,b,c)$ and $f' = (a',b',c')$ are \emph{equivalent} if there is an $M_{o} \in \psl$ so that $M_{o}\cdot f = f'$. This action leaves the discriminant $\Delta$ invariant, that is $\Delta(f) = \Delta(f')$ if $f$ and $f'$ are equivalent, and the equivalence class of a form $f$ is denoted by $[f]$. Therefore $\psl$ acts on the set of non-degenerate BQFs of the same discriminant, say a square-free integer $\Delta$, denoted by $\mathcal{F}(\Delta)$. In the search for a canonical representative of each class, in \cite{disquisitiones}, Gauss have defined a form $f = (a,b,c)$ to be \emph{reduced} whenever $\left|\sqrt{\Delta(f)} - 2\left|a\right|\right|<b<\sqrt{\Delta(f)}$. It turns out that whenever $f$ is definite (positive of negative) then there is a unique reduced BQF in its class, $[f]$. If $f$ is indefinite, then there are at least two reduced forms in $[f]$. Gauss have also given an algorithm which takes any non-degenerate BQF as an input and produces an equivalent reduced BQF, described in \cite{disquisitiones}.

One must note that there are other inequivalent notions of being reduced. For instance, Lagrange defines a BQF, say $f = (a,b,c)$, to be reduced when $|b|\leq a \leq c $ and $b \geq 0$ if either $a = c$ or $|b| = a$. Just as in the case of Gauss, every $\psl$-class of a positive definite BQF contains a unique Lagrange reduced form which can be found by the an algorithm that takes any non-degenerate BQF as input and produces an equivalent Lagrange reduced BQF, see the proof of \cite[Theorem~2.8]{cox/primes/of/the/form}. According to Zagier, an indefinite BQF is called reduced if both $\sqrt{\Delta(f)} < b < \sqrt{\Delta(f)} + 2a$ and $\sqrt{\Delta(f)} < b < \sqrt{\Delta(f)} + 2c$ are satisfied by $f = (a,b,c)$. We refer to \cite[\S~13]{zagier/zetafunktionen/quadratische/zahlkorper} for further details. 
\subsection{Ideal classes in quadratic number fields}
A field $K$ is called a \emph{number field} whenever it is a finite extension of $\QQ$. The dimension of $K$ as a vector space over $\QQ$ is called the degree of the extension and denoted by $[K:\QQ]$. $K$ is called \emph{quadratic} quadratic whenever $[K:\QQ] = 2$. In this case, one can always find a square-free integer $d$ with the property that $K = \QQ(\sqrt{d})$. If $d>0$, then we say that the extension is \emph{real} and if $d<0$ then the extension is called \emph{imaginary}. The map $\ol{\cdot}\colon K \lra K$ sending an element $\alpha = a + b\sqrt{d}$ to $\ol{\alpha} := a - b\sqrt{d}$ is the only non-trivial automorphism of $K$ which, together with identity, form the Galois group of $K$.

An element $\alpha \in \QQ(\sqrt{d})$ is called an \emph{algebraic integer} if it is a root of a monic polynomial $p_{\alpha}(x) \in \ZZ[x]$. It is not so hard to see that the set of algebraic integers in $\QQ(\sqrt{d})$, which will be denoted by $\oo_{d}$, depends on the square-free integer $d$ in the following manner: if $d$ is congruent to $1$ modulo $4$ then $\oo_{d} = \ZZ +\frac{1 + \sqrt{d}}{2} \ZZ$, and $\oo_{d} = \ZZ + \sqrt{d} \ZZ$ otherwise. In particular, $\oo_{d}$ is a Dedekind ring. The property of $\oo_{d}$ being a unique factorization domain (UFD) is equivalent to $\oo_{d}$ being a principal ideal domain(PID). And for any given ideal, $\mathfrak{a}$, of $\oo_{d}$ there are at most $2$ elements $\alpha$ and $\beta$ of $\oo_{d}$ so that the ideal generated by $\alpha$ and $\beta$, denoted $\langle \alpha,\beta \rangle$, is equal to $\mathfrak{a}$. This means fractional ideals of $K$, that is two dimensional $\ZZ$-modules, $\mathfrak{a}$, of $K$ for which there is an element $k_{\mathfrak{a}} \in \ZZ$ so that $k_{\mathfrak{a}} \cdot \mathfrak{a}$ is an ideal of $\oo_{d}$, can also be generated by at most $2$ elements. Given a fractional ideal $\mathfrak{a} = \langle \alpha,\beta \rangle$, we say that the pair $(\alpha,\beta)$ is \emph{oriented} if $({\ol{\alpha}\beta - \ol{\beta}\alpha })/{\sqrt{d}}>0$. The set of oriented fractional ideals is an abelian group under multiplication, denoted by $I^{+}(K)$. The subset of principal fractional ideals, that is subsets of the forms $\alpha \oo_{d} \subset K$ where $\alpha \in K$, is a subgroup of $I^{+}(K)$ denoted by $P(K)$, and the quotient $H^{+}(K) = I^{+}(K) / P(K)$ is called the \emph{narrow ideal class group} of $K$. The \emph{ideal class group} of $K$, denoted $H(K)$, is defined as the quotient $I(K) / P(K)$; where $I(K)$ is the multiplicative group of fractional ideals of $K$. The group $H(K)$ is naturally a subgroup of $H^{+}(K)$. $H(K)$ is of index $2$ if $\oo_{d}$ admits an element of norm $-1$ and $H(K)$ and $H^{+}(K)$ are isomorphic if there are no elements of norm $-1$ in $\oo_{d}$. 

\subsection{\c{C}arks and ideal classes}

\begin{figure}[H]
	\centering
	\includegraphics{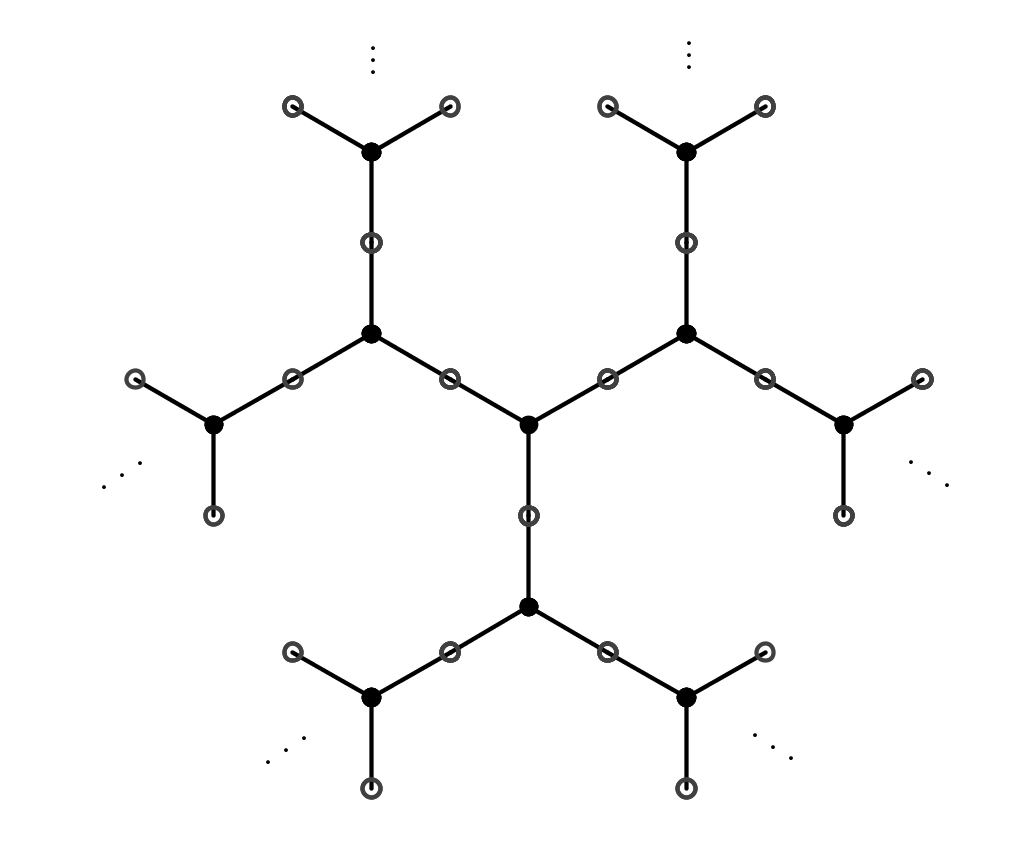}
	\caption{The Farey tree.}
	\label{fig:farey/tree}
\end{figure}

The orbit of the geodesic $\gamma$ in the upper half plane joining $\sqrt{-1}$, marked with $\circ$, to  $e^{2 \pi \sqrt{-1}/3}$, marked with $\bullet$, is a bipartite ribbon graph\footnote{A bipartite graph is a graph whose vertices can be decomposed into two disjoint sets so that no two edge belonging to the same set are joined by an edge. A ribbon graph is a graph together with a cyclic ordering of edges emanating from each vertex.} in $\HH$, where the vertices decompose naturally into two sets $V_{\circ} = \{\mbox{orbits of }\sqrt{-1} \}$ and $V_{\bullet} = \{\mbox{orbits of }e^{2 \pi \sqrt{-1}/3} \}$. The orientation is induced from the orientation of $\HH$. This graph is called the (bipartite) Farey tree and denoted by $\mathcal{FT}$, see Figure~\ref{fig:farey/tree}. Vertices of type $\circ$ are always of valency (or order or degree) $2$ and vertices of type $\bullet$ are always of valency $3$. 

%\begin{figure}[h!]

\begin{figure}[H]
	\centering
	\includegraphics{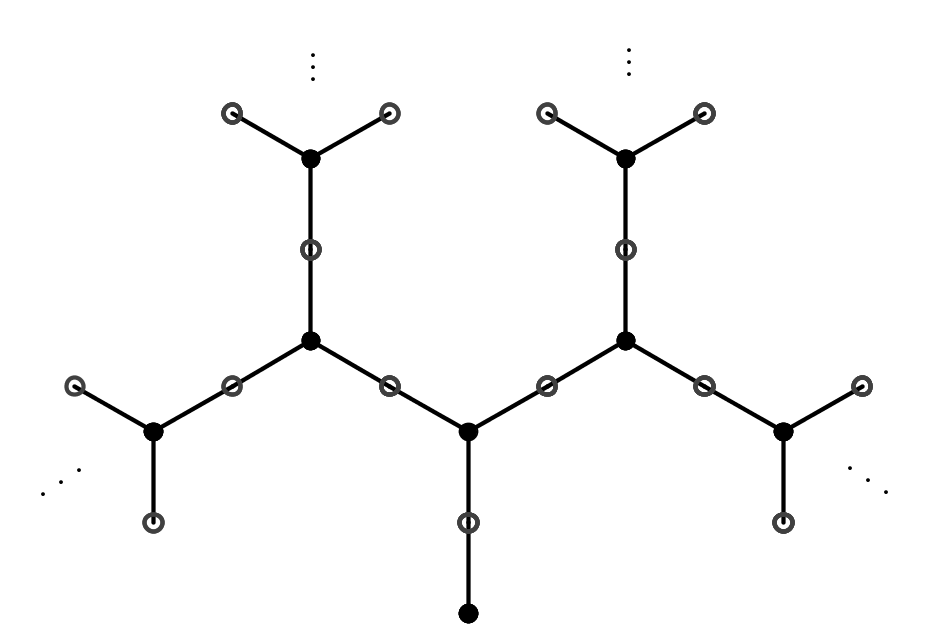}
	\qquad
	\includegraphics[scale=0.5]{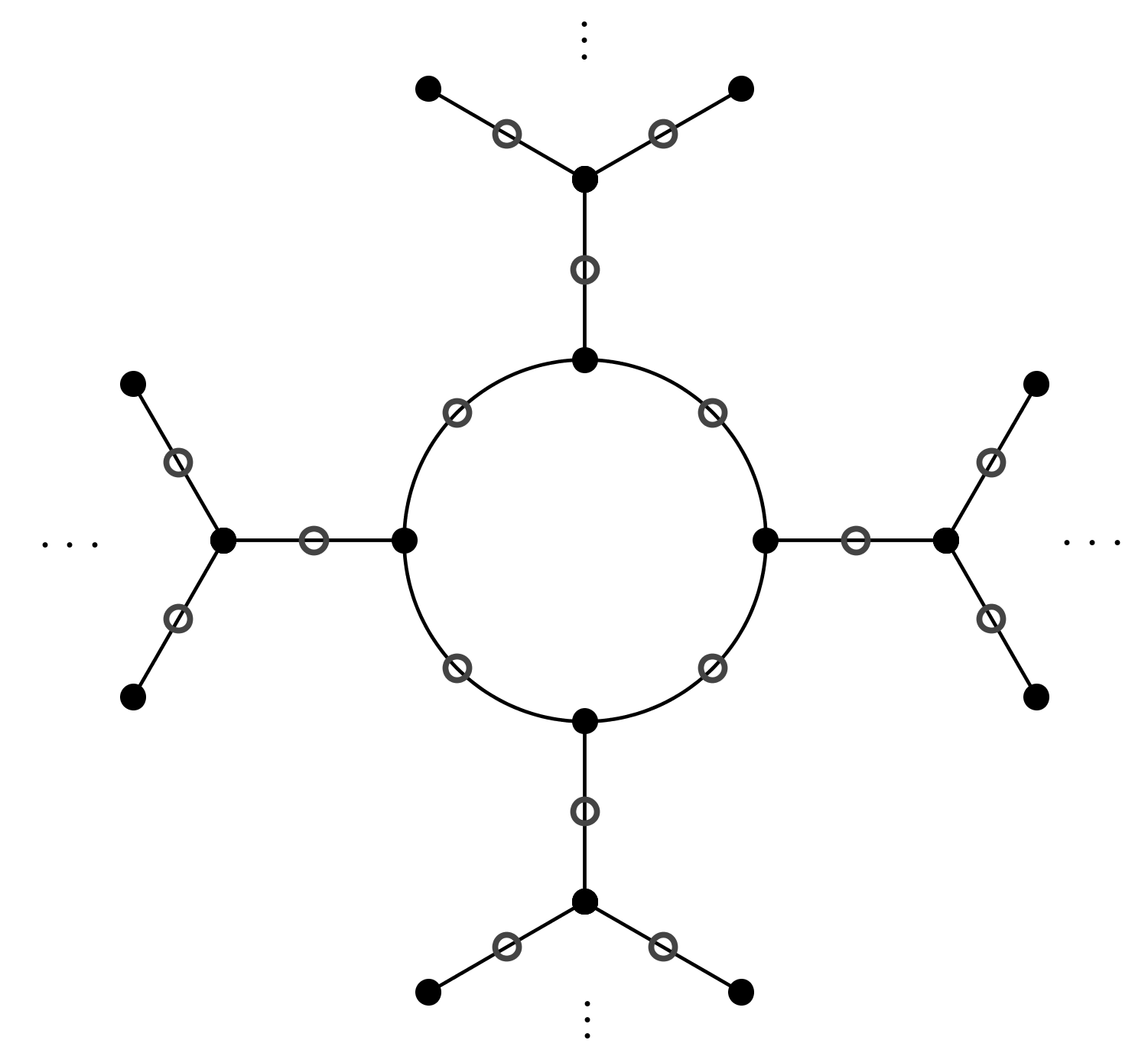}
	\qquad
	\includegraphics[scale=0.5]{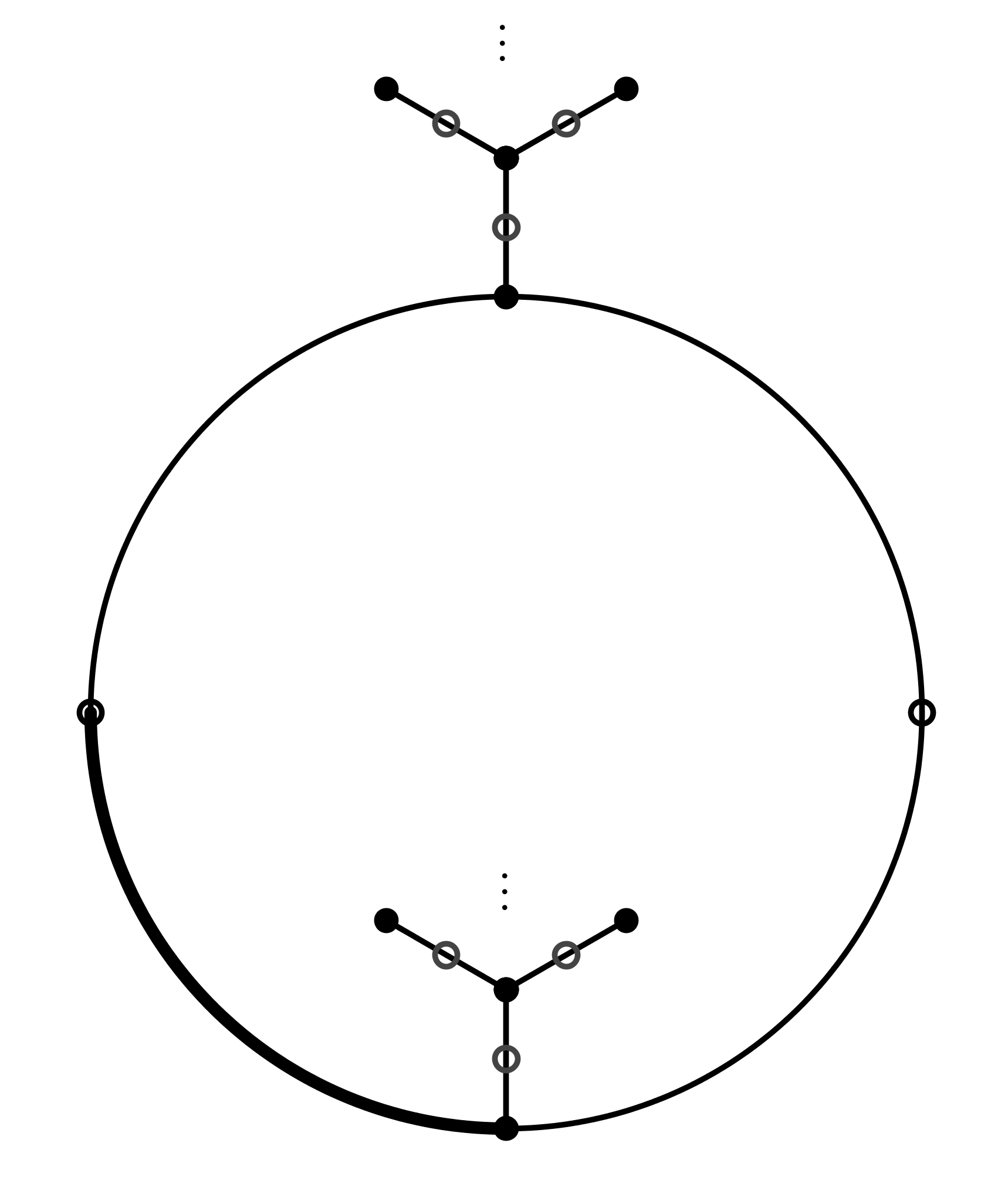}
	\caption{Elliptic, parabolic and hyperbolic \c{c}arks, respectively.}
	\label{fig:classification/of/carks}
\end{figure}

Given any element $W \in \psl$, by restricting the action of $\psl$ on $\HH$ to an action on $\mathcal{FT}$ we can consider the quotient of $\mathcal{FT}$ by the subgroup generated by $W$. The quotient is a bipartite ribbon graph which we refer to as a \emph{\c{c}ark}. The classification of elements of $\psl$ reflects itself in its \c{c}ark, see Figure~\ref{fig:classification/of/carks}.

%\begin{figure}[h!]

When $W$ is elliptic of order $2$ or $3$ we obtain a rooted tree. If $W$ is parabolic then the corresponding \c{c}ark has a unique cycle, named \emph{spine}, which is assumed to be oriented counter-clock-wise and finitely many rooted trees attached to this cycle each of which is a rooted tree attached to vertices of type $\bullet$ of the spine. Such components are called Farey components of the \c{c}ark. They all point outwards. When $W$ is hyperbolic, once again the graph has a unique cycle, named \emph{spine}, having finitely many vertices. To each vertex of type $\bullet$ there is a Farey component attached. In this case, Farey components must point both inwards and outwards. 

As a result of the construction the cosets $\gamma \cdot \langle W\rangle \in \psl / \langle W \rangle$ correspond in a one to one fashion to edges of the \c{c}ark $\mathcal{FT}/\langle W \rangle$. So we conclude that there is a one to one correspondence between the set of edges of the {\cark} corresponding to the form $f_{W}$ and element of the equivalence class $[f]$, see \cite[\S~2]{UZD} for details.

\section{Software Design}

\subsection{Components}

InfoMod consists of a shared library coded in C for computationaly heavy functions, some convenience wrapper classes for Matlab and Python to facilitate interfacing the library, and a visualization application written in ActionScript version 3.

The native code shared library implements the computation of the two methods of reduction (\c{c}ark and Gauss), enumeration of quadratic forms residing on a spine and computation of the signature of a spine. On the other hand wrapper classes to access the library from Matlab and Python environments provide an easy way to invoke library functions, while providing an Object Oriented representation for the involved quadratic forms. Source code for Matlab and Python libraries are available online \url{https://github.com/hayral/infomod}. Many functions which aren't computationaly demanding, such as evaluating a form at a point, obtaining the matrix representation of the form, or querying the spinality, primitivity and reducedness, are implemented in the wrapper classes, as the overhead of function invocation from dynamic library outweights the benefits for those cases.

\begin{figure}[H]
	\centerline{\includegraphics[scale=0.45]{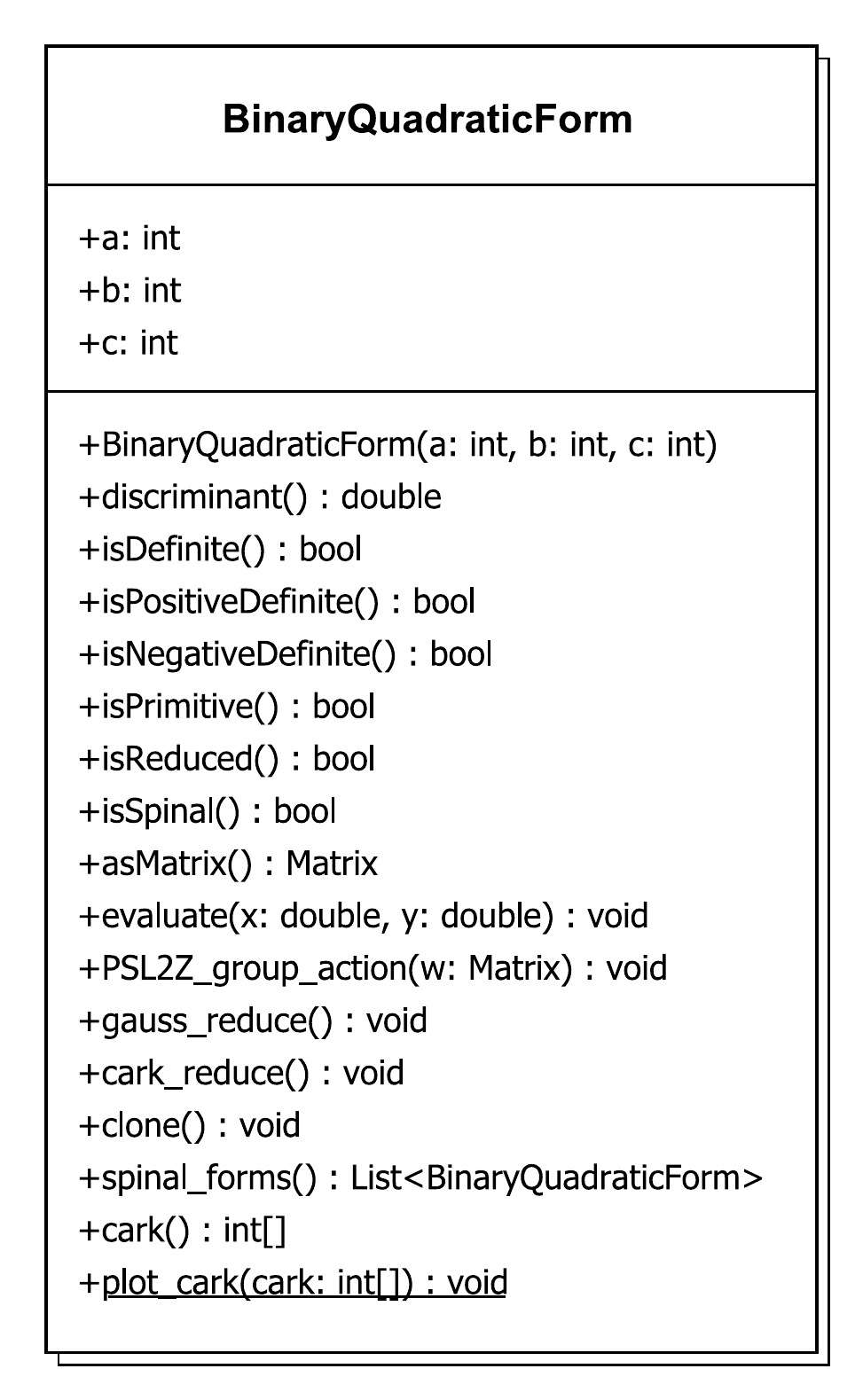}}
	\caption{BinaryQuadraticForm class in Matlab and Python wrapping low-level reduction library and providing some convenience functions. (parameter and return types seen on figure are for UML diagramming purposes, actual implementation types differ according to Matlab/Python environment )}
	\label{fig:bqfclass}
\end{figure}

\subsection{Visualization of the modular group and Sunburst}

While Matlab and Python libraries are designed to be used as part of other computational code, the interactive visualization application of InfoMod is designed as an exploration tool to provide insight. It uses the same code base as the libraries, but ported to ActionScript programming language; as we planned to deploy this application to web and mobile, invoking our native code shared library or relying on Matlab/Python was not an option. Before elaborating on the visualization aspect of the software, let us define the underlying mathematical structure that visualization is designed to explore.

Modular group is generated by two torsion elements $S = \SS$ and $L = \LL$ of orders $2$ and $3$, respectively. There are no relations among $S$ and $L$ hence we have the isomorphism $\psl \cong \ZZ/2\ZZ \ast \ZZ/3\ZZ$. So every element in $\psl$ can be written as a word in $S$, $L$ and $L^{2}$. Without loss of generality, we assume that the length of words, denoted by $\ell(W)$, are of minimal length\footnote{The length of a word is defined to be the number of letters ($L$, $L^{2}$ and $S$) appearing in the word.}. Given two words $W$ and $W'$ in $\psl$, by $W\cap W'$ we denote the word which is equal to the common initial part of the words $W$ and $W'$. For instance, for $W = (LS)^{2}(L^{2}S)^{3}LSL$ and $W' = (LS)^{2}(L^{2}S)^{3}L^{2}SLSL^{2}$ we have $W \cap W' = (LS)^{2}(L^{2}S)^{3}$.

\begin{figure}[H]
	\centerline{\includegraphics[scale=0.35]{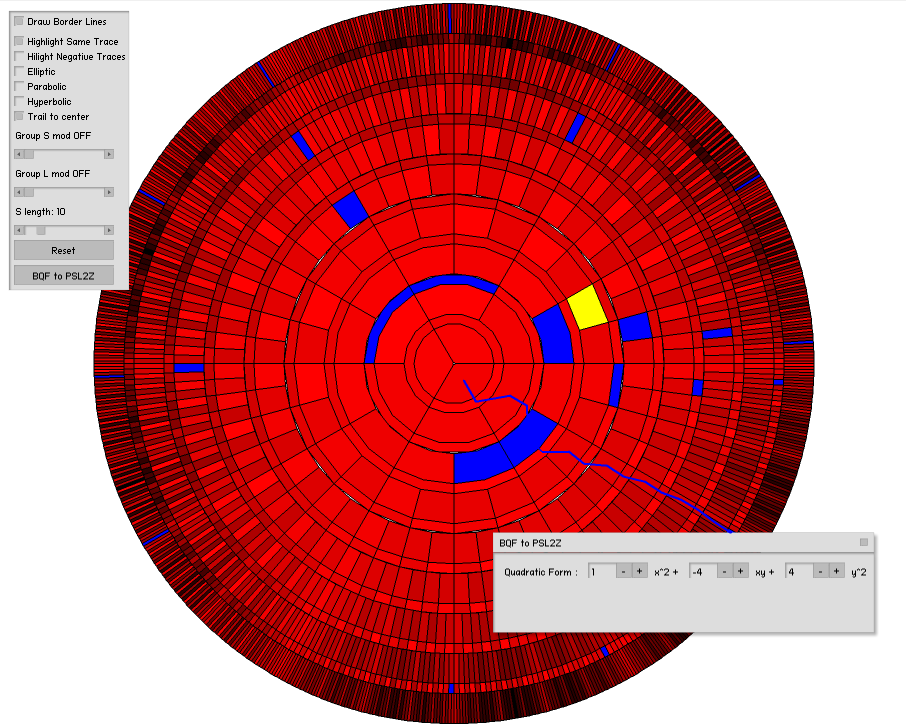}}
	\caption{Main display of InfoMod}
	\label{fig:infomod}
\end{figure}

Let us now construct the slit disk, which we denote by $\dd$. We start with a standard disk, split into three equal area pieces by lines from the center of the disk. The so obtained three cells are labeled as $I$, $L$ and $L^{2}$. To the circumference of this disk we glue an annulus which is also divided in $3$ pieces by prolonging the lines that were used to divide the initial disk. This time, the 3 new cells are identified with $S$, $LS$ and $L^{2}S$, respectively. The next step is to attach two annuli each of which is divided evenly into 6 pieces via first prolonging the existing lines and then adding 3 more lines so that the cell labeled $S$ is neighbors with the cells $SL$ and $SL^{2}$, and these are in turn neighbors with $SLS$ and $SL^{2}S$ in the second annulus. Similarly, the cell labeled $LS$ is neighbors with the cells $LSL$ and $LSL^{2}$, and these are in turn neighbors with $LSLS$ and $LSL^{2}S$ in the second annulus and finally the cell labeled $L^{2}S$ is neighbors with the cells $L^{2}SL$ and $L^{2}SL^{2}$, and these are in turn neighbors with $L^{2}SLS$ and $L^{2}SL^{2}S$ in the second annulus. Inductively, we obtain a disk, to which we will refer as the slit disk, where each cell is labeled with a word in $S$, $L$ and $L^{2}$ see Figure~\ref{fig:infomod}. As a result of the construction given above we have a one to one correspondence between elements of $\psl$ and cells in  $\dd$. 

If the a cell in $\dd$ is labeled as $W$ and ending with $L$ or $L^{2}$ then there are 4 cells surrounding it: three of them are labeled with $WL$, $WL^{2}$ and $WS$. Remark that $\ell(W \cap W') \geq \ell(W) - 1$ where $W'$ is one of $WL$, $WL^{2}$ or $WL^{2}$. The fourth neighbor, say $W''$, of $W$ satisfies $\ell(W \cap W'') \leq \ell(W)-2$ with $\ell(W) = \ell(W'')$. Similarly, if $W$ is a word ending with $S$, then there are 5 cells surrounding this cell. Three of them, which satisfies $\ell(W \cap W') \geq \ell(W)-1$, are labeled with $WL$, $WL^{2}$ and $WS$. The remaining two words $W''$ and $W'''$ satisfy $\ell(W) = \ell(W'') = \ell(W''')$ and $\ell(W \cap W'')<\ell(W)-2$ and $\ell(W \cap W''')<\ell(W)-2$.

% \subsection{Modular Group Sunburst}

% % % The main visualization of the application is a sunburst type diagram of the tree topology induced by the generators of $PSL_{2}(\mathbb{Z})$. 
% % % A disc, at the center of many other concentric layers, is divided to three equal sections representing the identity element of modular group and the first two powers of the third degree generator of the group, $L$ and $L^2$, which are respectively equal to $\bigl(\begin{smallmatrix}
% % % 1&-1 \\ 1&0
% % % \end{smallmatrix} \bigr)$ and $\bigl(\begin{smallmatrix}
% % % 0&-1 \\ 1&-1
% % % \end{smallmatrix} \bigr)$ in matrix form. A thin annulus encircles the inner most disc, and consists of three equal segments which are the elements of modular group corresponding to the inner ones multiplied with the degree 2 generator $S$ of the group, which is equal to $\bigl(\begin{smallmatrix}
% % % 0&-1 \\ 1&0
% % % \end{smallmatrix} \bigr)$ in matrix form. As the radius increase from here to the outermost boundary of the big disc composed of all the annulus segments, every annulus segment which is the result of a multiplication by the generator $S$ is followed by two other annulus segments, dividing the angular range to two equal parts, and corresponding to further multiplication by $R$ and $R^2$.

The main visualization of the application is the above described slit disc (sunburst) type diagram. When the user moves the mouse pointer over one of the cells on $\dd$, a blue line indicates the path to the center of the disc as it passes by the centers of the parent cells. Due to the exponential growth of the number of child nodes on a binary tree, this sunburst type visualization can only depict eight levels of depth for common mobile device and desktop screen resolutions.  Once a cell representing an element of $\psl$ is left clicked, a radial menu appears with options. Apart from the elements in the menu related to the visualization of the corresponding binary quadratic form along with its neighbors the button ``Move to Center'' is used to interact with elements of the modular group which are at depth nine or further. Namely, on click by conjugation, the chosen element becomes the center of the displayed $\dd$ and all the other elements are translated accordingly. In this fashion, one may travel \emph{arbitrarily} far from the actual center. 

%\note{hakan: bu yazdigim ne kadar mantikli bilmiyorum. Uygun dile dokmek lazim.}

\begin{figure}[H]
  \centering
  \subfloat{\includegraphics[scale=0.5]{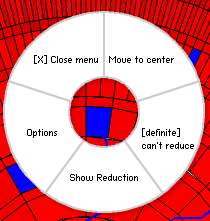}\label{fig:f1}}
  \qquad  \qquad
  \subfloat{\includegraphics[scale=0.3]{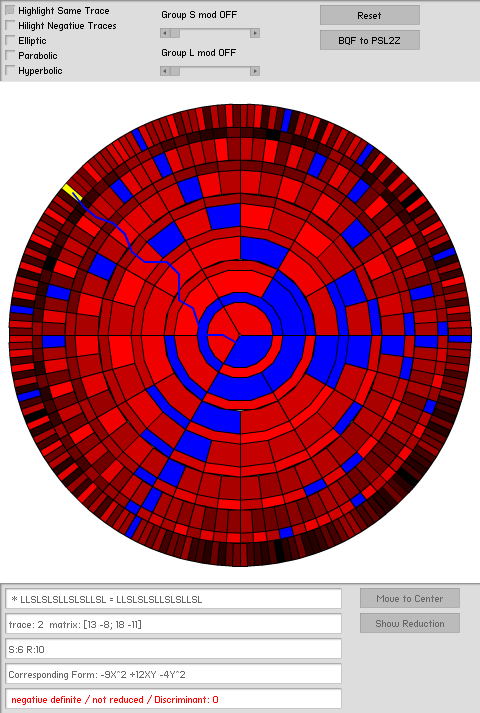}\label{fig:f2}}
  \caption{Radial menu on web application and main screen on mobile application}
\end{figure}

 \paragraph{Challenges implementing the web and mobile visualization app} ActionScript isn't designed with scientific computing or high performance computing as primary concern, but it has all the features of a modern managed object oriented language \cite{crawford2006actionscript}, and provides access to hardware accelerated graphics capabilities of the underlying operating system. %The other two constituents of the platform are the Flash runtime and the Flash SDK which includes the compiler and some libraries. 
Flash runtime is generally deployed in the form of a browser plugin (Flash Player), but standalone applications can be deployed as  independent executables without requiring being embedded in a browser. %Runtime provides the virtual processor for the Flash bytecode just like what JRE does for Java, and establishes the interface with the operating system for user interaction, network and file system access. Flash SDK consists of some documentation and samples but most importantly the compiler which converts the ActionScript source code to Flash bytecode binary, and the library files which provide callable handles to the services and functionalities provided by the runtime.

ActionScript and the Flash Runtime has no direct access to mathematical computing libraries; therefore one of the challenges was the necessity to implement all mathematical structures and computation codes from scratch; including a class to represent a binary quadratic form, a class to represent an element from the modular group, and minimal linear algebra to operate on the matrix representation of elements of modular group\footnote{Standard ActionScript libraries do provide a native $3 \times 3$ matrix class, but in our case we needed some extra functions like trace,transpose and multiplication with specific constants (i.e. generators of modular group and some of their compositions) which reduces to simple swap and sign change of values when optimized.}.

Another challenge was the code being run by a JIT compiling VM with automatic garbage collection which is the Flash Runtime. Even though we didn't need C/C++ like native code speed, having numerous visual objects (up to the order of $2^{12}$) visible on the screen with many destroyed and some newly created at each user interaction cause a heavy burden on the garbage collector; therefore we employed some object pooling and caching methods to maintain interactive speeds and minimize stutter caused by garbage collector.
Another challenge was to port the same application to mobile platforms, namely Android and iOS. AIR compiler provides platform independence to some degree, but user experience is not the same due to differences on hardware features. One of the two main manifestations of this are the resolution and aspect ratio differences of different mobile devices requiring the size and placement of visual objects to be flexible and independent of screen size (including fonts and text rendering which is not as fast and flexible on mobile platforms); and the other being the relatively limited resources in terms of computing power and memory.

% % % For every annulus segment representing an element of modular group which can be represented in matrix form as $\bigl(\begin{smallmatrix}
% % % a&b \\ c&d
% % % \end{smallmatrix} \bigr)$, there is a corresponding binary quadratic form $AX^2+BXY+CY^2$ which can be obtained by the relations $A= \frac{c}{\delta}, B= \frac{d-a}{\delta}, C= \frac{-b}{\delta}$. Furthermore classes of indefinite binary quadratic forms have a neighborhood in form of a planar graph topology, where the reduced forms form a central cycle (called spine in \c{c}ark terminology), and some tree structures of infinite depth emerge from the elements of the cycle either inwards or outwards the cycle \cite{reduction}.

\subsection{Visualization of the \c{c}arks and geodesics}
Let us now fix a hyperbolic element $W \in \psl$. The \c{c}ark visualization consists of edges forming the spine (which are called \emph{semi-reduced}) and the rooted Farey trees  attached to it. The correspondence between edges of $\FT/\langle W \rangle$ and forms $f$ that are equivalent to $f_{W}$ can be investigated by hovering the mouse cursor over, thus highlighting the edge green and updating the top left corner green label with the corresponding forms expression.  The form corresponding to $W$ is indicated as a red edge, and has a label for its expression on top left.

%\begin{figure}[H]
%\centerline{\includegraphics[scale=0.35]{2.png}}
%\caption{\c{C}ark visualization}
%\end{figure}

%\begin{figure}[H]
%\centering
% --------------------------------------------

%\begin{subfigure}{.5\textwidth}
%  \centering
%  \includegraphics[height=.5\linewidth]{2.png}
%  \caption{As rendered by interactive visualization application}
%\end{subfigure}
%
%\begin{subfigure}{.5\textwidth}
%  \centering
%  \includegraphics[height=.5\linewidth]{cark-14+2+1.pdf}
%  \caption{As rendered by Matlab wrapper}
%
%\end{subfigure}%
%
%
%\caption{Visualizations of \c{c}ark for form (-14,2,1)}
%\end{figure}

% --------------------------------------------

\begin{figure}[H]
  \centering
  \subfloat{\includegraphics[height=.5\linewidth]{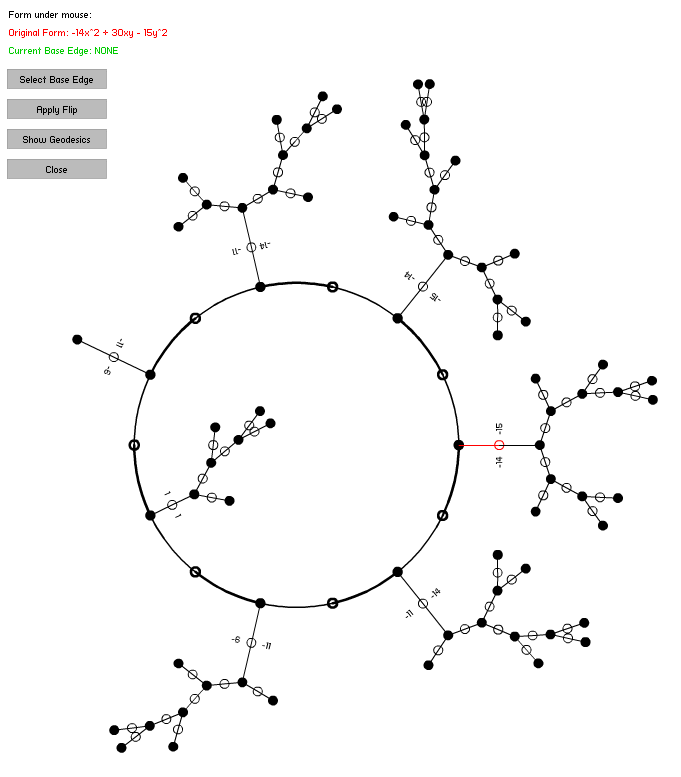}} 
  \qquad  \qquad
  \subfloat{\includegraphics[height=.45\linewidth]{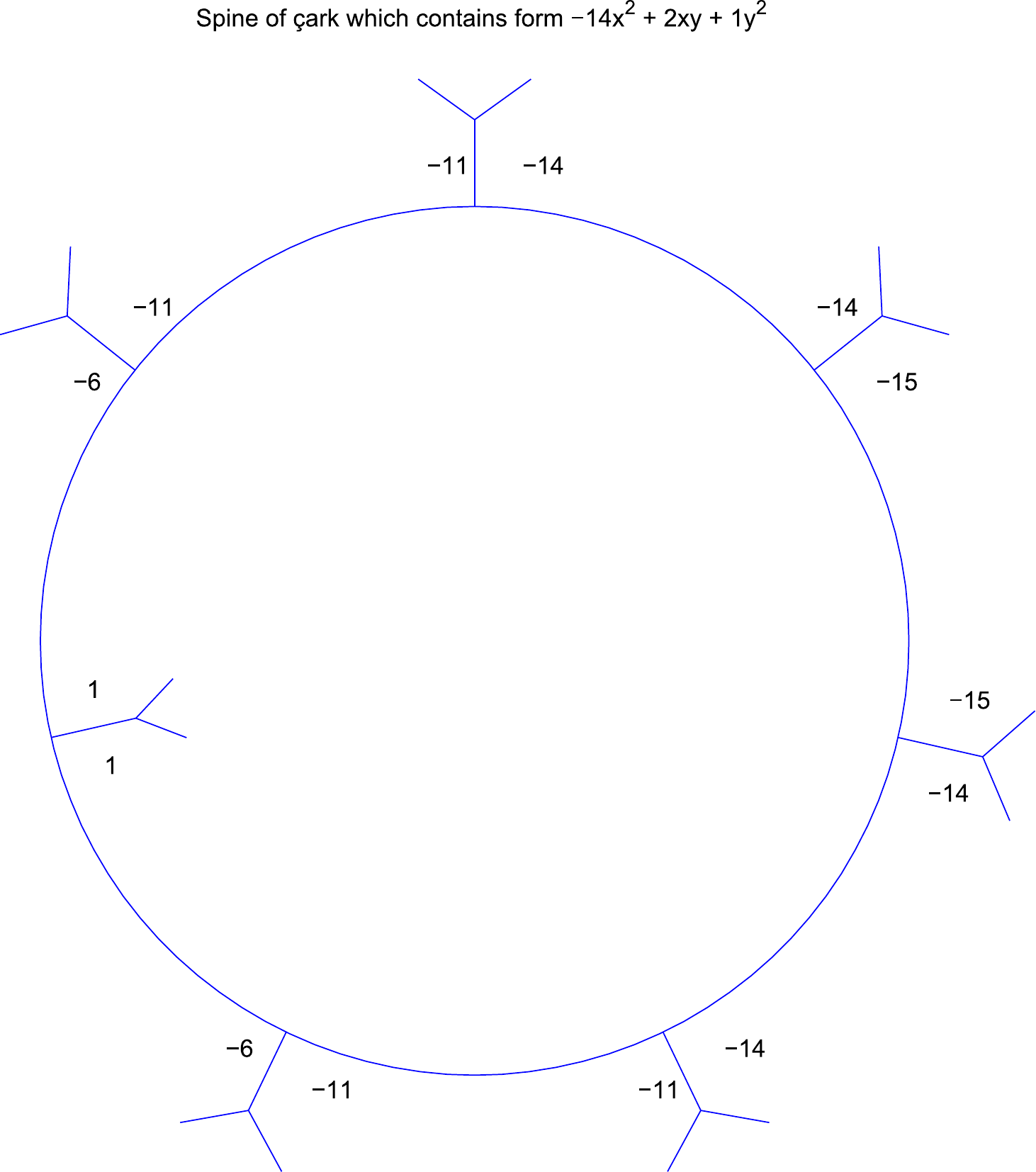}} 
  \caption{Visualizations of \c{c}ark for form (-14,2,1). (left) As rendered by visualization application, and some branches expanded through interaction (right) As rendered by plot\_cark() method of Matlab wrapper class}
\end{figure}

Every non-reduced form in the equivalence class (i.e. edges not on the spine) can be clicked on to expand the child edges, and the whole graph re-scales and positions itself to fit on available screen space; yet as the trees are infinite deeper child edges get smaller angular ranges to position themselves.

Every quadratic form on the spine of the \c{c}ark also corresponds to a geodesic in upper half-plane. Once an edge is selected on the \c{c}ark visualization, the geodesic of the form corresponding to selected edge can be plotted, both in $\HH$ or in the unit disk, $\DD = \{z\in \CC \colon |z|<1\}$.

%\begin{figure}[H]
%\centerline{\includegraphics[scale=0.3]{4.png}}
%\caption{Visualization of the geodesics}
%\end{figure}

\begin{figure}[H]
  \centering
  \subfloat{\includegraphics[height=.35\linewidth]{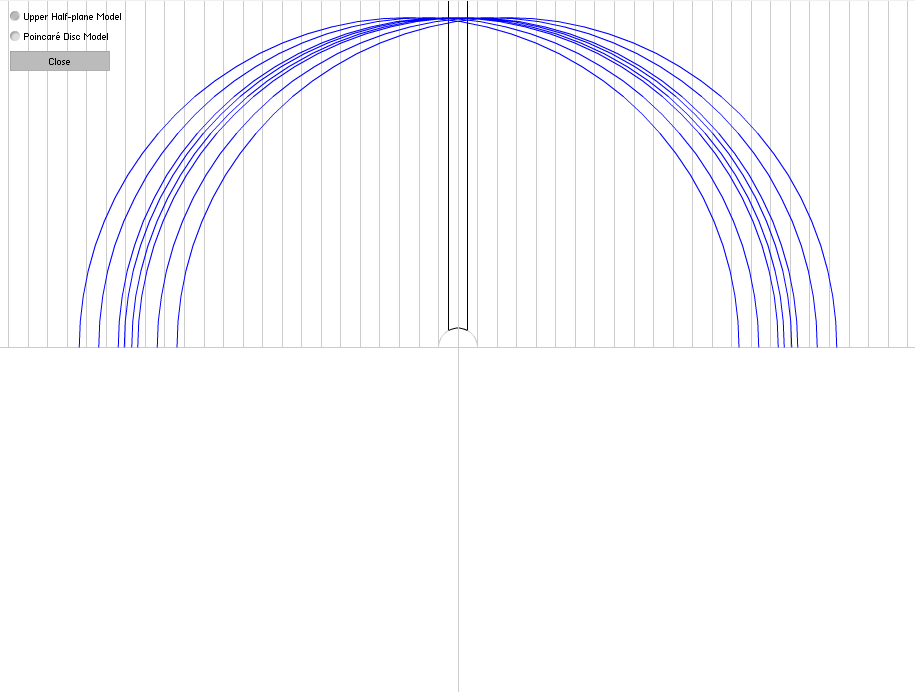}} 
  \qquad  \qquad
  \subfloat{\includegraphics[height=.35\linewidth]{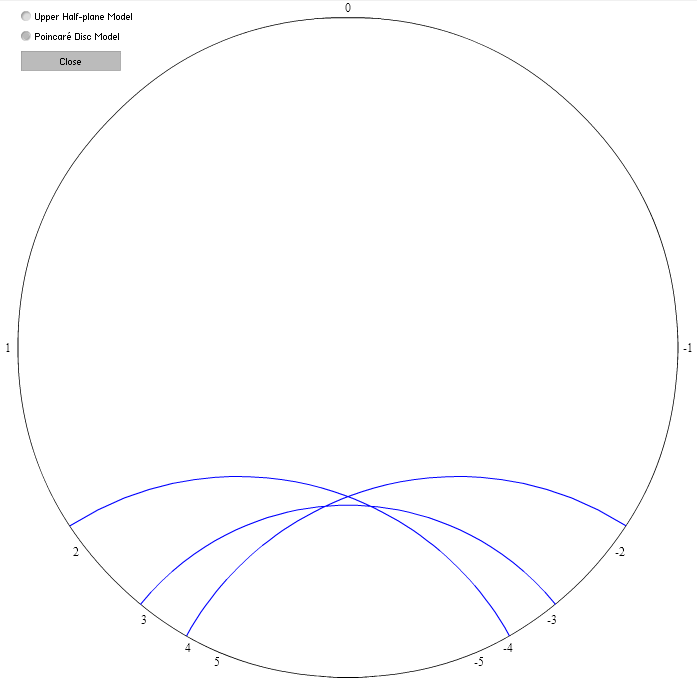}} 
  \caption{Visualization of the geodesics (left) Upper half plane model (right) Disc model}
\end{figure}

% % % \subsection{Intervals covered by a base edge}

\begin{figure}[H]
  \centering
  \subfloat{\includegraphics[width=.45\linewidth]{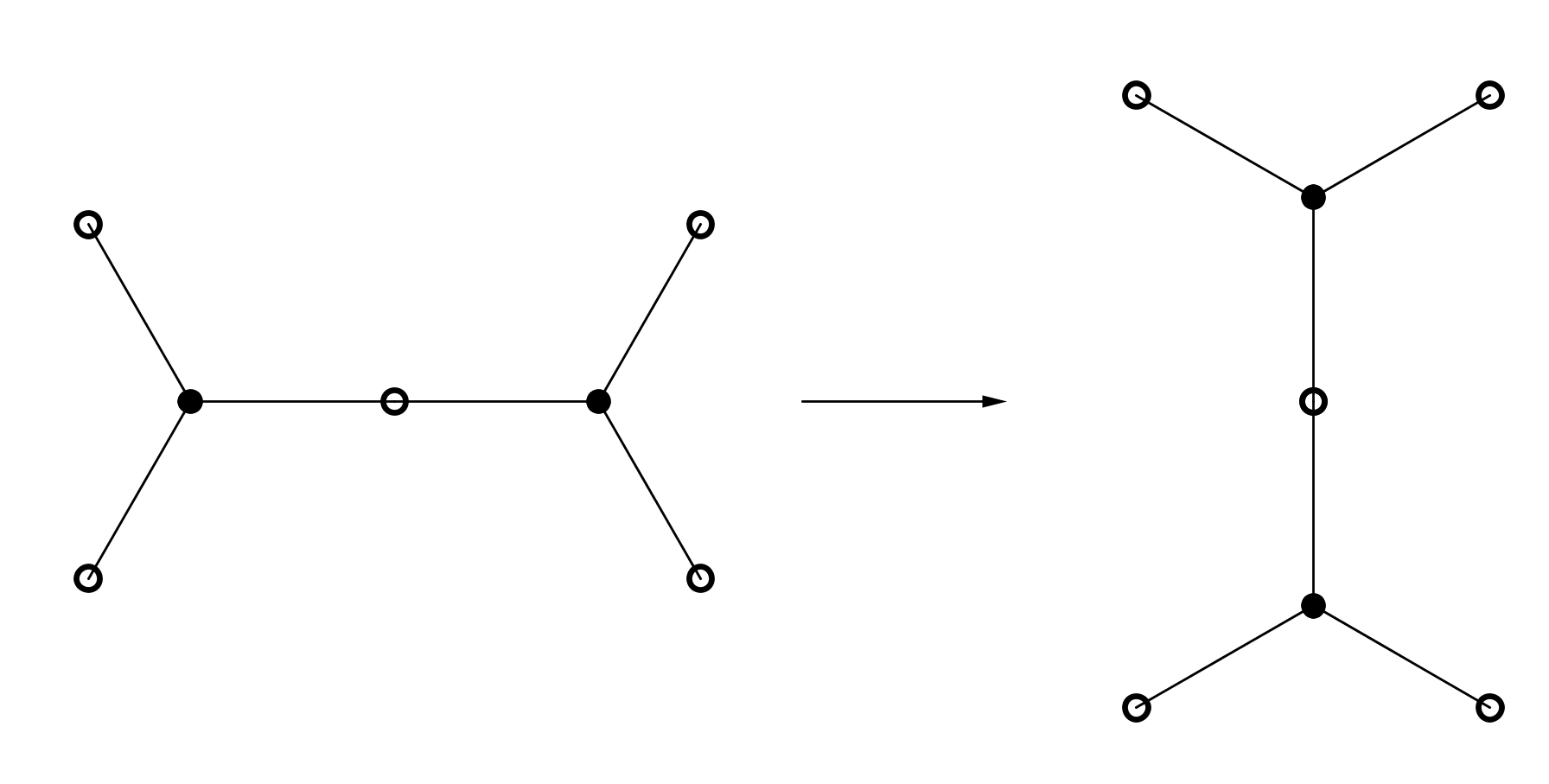}} 
  \qquad  \qquad
  \subfloat{\includegraphics[width=.40\linewidth]{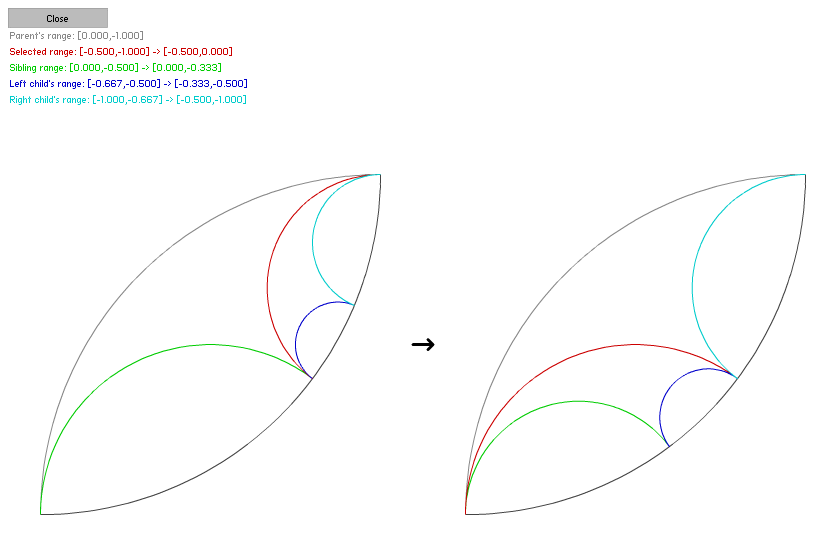}} 
  \caption{(left) Flip action on ribbon graphs (right) Visualization of intervals exchanged on the boundary of the disk for a specific flip}
   \label{fig:flip}
\end{figure}

\paragraph{Flips} There is a close relationship between mapping class groups of surfaces of genus $g$ with $n$ punctures and the groupoid whose objects are bipartite ribbon graphs of genus $g$ with $n$ punctures and morphisms are flips (or Whitehead moves or HI moves), see Figure~\ref{fig:flip}. Indeed, flips do not change the invariants $g$ and $n$ and act transitively on the set of ribbon graphs of fixed $g$ and $n$. Then, as a result of a generalization of Dehn-Nielsen-Baer theorem, \cite[Theorem~4.8]{farb/margalit/MCG},   we obtain the required injection. InfoMod can make explicit computations using flips in the case of annulus.

%\begin{figure}[h!]
%	\centering
%	\includegraphics{flip.png}
%	\caption{Flip action on ribbon graphs}
%	\label{fig:flip}
%\end{figure}

\section{Representation problem}

Representation problem deals with the existence and values of integer solutions, that makes a given binary quadratic form equal to a chosen integer. For the positive definite forms if integer solutions exist, they are finitely many; on the other hand for indefinite case if an integer solution exists, then there is a family of countably infinite pairs of integers. Each solution pair $(x,y)$ can be obtained by multiplying a non-trivial solution with an automorphism of the form of appropriate order. 

It must be mentioned that there is a non-exact algorithm presented on the web-page ``https://www.alpertron.com.ar/quad.java''. In addition, in \cite{lagrange} Lagrange has put forward a method to solve similar equations which run in exponential time, which are then further optimized in \cite{sawilla/silvester/williams}. Some related algorithms have been investigated in \cite{lagarias/worst/case/complexity/bound}. Pell equation is another class of such equations whose theory is well developed, \cite{jacobson/williams/solving/pell/equation}.  In \cite{jacobson/williams} equations of the form $aX^{2} - bY^{2} = N$ are discussed. Algorithms concerning the case $N = -1$ (sometimes referred to as non-Pellian equation), have been discussed by Lagarias in \cite{lagarias/computational/complexity/of/pell}.

In \cite{lagarias/computational/complexity/of/pell} Lagarias proposes a modification to classical reduction method, through which he established a worst case time complexity of $\mathcal{O}(n\mu(n))$ with $\mu(n)$ being the time complexity of the chosen n-bit multiplication algorithm, yielding $\mathcal{O}(n^3)$ in case regular multiplication is employed. In \cite{Buchmann} Buchmann improves the complexity analysis of Lagarias to $\mathcal{O}(n^2)$ without changing the algorithm, using the observation that if for a form $f$ neither $f,\rho(f)$ or $\rho^2(f)$ \footnote{with $\rho$ being the reduction operator, defined as $fU(f)$ and $U(f)=\begin{pmatrix}0 & -1 \\ 1 & s(f) \end{pmatrix}$} is reduced then the absolute value of the first coefficient ($a$) must be at least twice the absolute value of the last coefficient ($c$). Eventhough a faster reduction algorithm of $\mathcal{O}(log(n)\mu(n))$ time complexity is presented earlier by Schönhage \cite{schonhage}, the asymptotic acceleration surpasses the other algorithms only when the number of binary digits of the coefficients of the form to be reduced is in the order of $10^4$ or greater.

The representation problem is closely related to class group computations in quadratic number fields. For instance, in \cite{shanks/class/number}, Shanks developed a new method for the computation of class numbers of quadratic number fields, which further resulted in a new factorization of large integers. The thesis of Jacobson, \cite{jacobson/phd/thesis}, discusses such algorithms in detail.  Similar results for general number fields are presented in \cite{cohen/diaz/olivier}.

Here, we present an algorithm to solve the representation problem given a binary quadratic form and an integer to solve for. The algorithm relies on two things: first to find another form which is part of a face with label equal to given integer if such a face exists (or to show that such a face doesn't exist), and second to compute the path leading from this form found to the form given.

% Algorithm
%\begin{algorithm}[t]
%\SetAlgoNoLine
%\KwIn{A binary quadratic form $F$ and the integer $N$ for which to solve $F$}
%\KwOut{An integer pair $X,Y$ if there is solution, or $NULL$ if there isn't}
%
%
%$index$ = 0; $FreNum_{\alpha}$ = -1\;
%\Repeat{$FreNum_{\alpha} > -1$}{
%        $Rnd_{\alpha}$ = Random($ID_{\alpha}$, $index$)\;
%        $Found$ = $TRUE$\;
%        \For{each node $\beta$ in $\alpha$'s two communication hops
%    }{
%      $Rnd_{\beta}$ = Random($ID_{\beta}$, $index$)\;
%      \If{($Rnd_{\alpha} < Rnd_{\beta}$) \text{or} ($Rnd_{\alpha}$ ==
%          $Rnd_{\beta}$ \text{and} $ID_{\alpha} < ID_{\beta}$)\;
%      }{
%        $Found$ = $FALSE$; break\;
%      }
%        }
%     \eIf{$Found$}{
%           $FreNum_{\alpha}$ = $index$\;
%         }{
%           $index$ ++\;
%     }
%      }
%\caption{Solver for Representation Problem}
%\label{alg:one}
%\end{algorithm}

\begin{algorithm}
\caption{Solver for Representation Problem}
\label{alg:algo1}
\begin{algorithmic}[1] 
\Require $f_{start}$ is form object, $N$ is an integer

\Function{SolveForm}{$f_{start}$,$N$}

\State $path_1 \gets$ \Call{CarkReducePath}{$f_{start}$} \Comment{path from $f_{start}$ to spine}
\State $entryPoint_1 \gets path_1($ \Call{len}{$path_1$} $)$ \Comment{last form on $path_1$ (it's on spine)}
\State $spineForms \gets$ \Call{RevolveAroundSpine}{$entryPoint$} \Comment{all forms on spine}
\State $formsToSearch \gets \emptyset  $
\State $ L \gets \begin{pmatrix} 1&-1\\ 1&0 \end{pmatrix} , S \gets \begin{pmatrix} 0&-1\\ 1&0 \end{pmatrix}  $ \Comment{$\psl$ generators}
\State $e_1 \gets (1,0) , e_2 \gets (0,1)$ \Comment $e_1,e_2$ : standard oriented basis of $\ZZ^2$

\For{each $f$ in $spineForms$} \Comment neighbours of spine
\State $f_1 \gets f \cdot L$ \Comment $\psl$ group action
\State $f_2 \gets f \cdot L^2$ 

\If{$(f_1(e_1) \times f_1(e_2) > 0) \land (f_1(e_1) \times N > 0)$ }
%	\Call{Append}{$formsToSearch,f_1$}
\State $formsToSearch = formsToSearch \cup \{f_1\}$
\EndIf

\If{$(f_2(e_1) \times f_2(e_2) > 0) \land (f_2(e_1) \times N > 0)$ }
%	\Call{Append}{$formsToSearch,f_2$}
\State $formsToSearch = formsToSearch \cup \{f_2\}$
\EndIf

\EndFor

\State $temp \gets \emptyset$ 

\While{$formsToSearch \neq \emptyset$} \Comment breadth-first search
\For{each $f \in formsToSearch$}
\If{$ ( f(e_1) = N) \lor (f(e_2) = N) $ }
	\State $f_{found} \gets f$ \Comment $f_{found}$: first form on the face with label $N$
	\State \textbf{break} \Comment leave the while loop
\Else

\State $f_1 \gets f \cdot S \cdot L$ \Comment $\psl$ group action
\State $f_2 \gets f \cdot S \cdot L^2$ 

\If{$(abs(f_1(e_1)) < abs(N)) \land (abs(f_1(e_2)) < abs(N))$ }
	%\Call{Append}{$temp,f_1$}
	\State $temp =  temp \cup \{f_1\}$
\EndIf
\If{$(abs(f_2(e_1)) < abs(N)) \land (abs(f_2(e_2)) < abs(N))$ }
	%\Call{Append}{$temp,f_2$}
	\State $temp =  temp \cup \{f_2\}$
\EndIf

\EndIf

\EndFor
\State $formsToSearch \gets temp , temp \gets \emptyset$ \Comment swap the lists

\EndWhile

\If{$targetForm = \emptyset$}
	\State \textbf{return} $\emptyset$ \Comment no integer solution
\Else
	\State $path_2 \gets$ \Call{CarkReducePath}{$f_{found}$} \Comment{path from $f_{found}$ to spine}
	\State $entryPoint_2 \gets path_2($ \Call{len}{$path_2$} $)$ \Comment{last form on $path_2$ (it's on spine)}
	\State $path_3 \gets$ \Call{PathOnSpine}{$entryPoint_2,entryPoint_1$} \Comment{path on spine between entry points} %from $entryPoint_2$ to $entryPoint_1$}
	\State $path_4 \gets path_2 + path_3 + $  \Call{ReversePath}{$path_1$} \Comment{path from $f_{found}$ to $f_{start}$}
	\State $M \gets$ \Call{PathToMatrix}{$path_4$}
\If{$f_{start}(e_1 \cdot M) = N$}
\State \textbf{return} $e_1 \cdot M$
\Else
\State \textbf{return} $e_2 \cdot M$
\EndIf

\EndIf

\EndFunction
\end{algorithmic}

\end{algorithm}

Equipped with reduction and graph topology of \c{c}arks, the computation of integer solutions proceeds as follows: First the software performs the reduction while recording the path down to the spine of its \c{c}ark as a sequence of $\psl$ generators (algorithm 1,line 2), then it travels along the spine and this time it records the forms belonging to spine as a sequence (l.4). Depending on the sign of the chosen integer to solve the quadratic form for, software generates all non-spinal neighbours of recorded spinal forms, but records only those with coefficients having same sign as the chosen integer (l.8-17). This list of non-spinal neighbours of the same sign with desired integer, serves as the root nodes of a simultaneous breadth-first search along the graph of the spine (l. 18-36). This breadth-first search travels away from the spine (inwards or outwards depending on the sign). Even though every form yields two new forms once expanded on the direction of increasing distance, the number of accumulated forms does not increase exponentially; the reason for this boils down to the coefficients of neighbouring forms growing or decreasing (depending on the travel being directed towards or away from the spine) at a non-linear speed when multiple forms are traversed. At each iteration of breadth-first search software compares the absolute value of the coefficients of the newly expanded forms with the absolute value of the integer being searched, and discards the form if the first or last coefficient of the form is greater than this in absolute value (l.27-32). If all the forms on breadth-fist search list gets discarded, the loop terminates (l.19) and algorithm doesn't return an integer pair, indicating that there are no integer solutions (l.37-38). If the search finds a form adjacent to a face with chosen label, it breaks from the loop after marking the found form (l.21-23), and the path from the found form and its entry point to spine is calculated (l.40-41). At this point as we have paths from both starting and found form to the spine; to build a continuous path from found form to starting form, we only need the path that connects the entry points of those paths to the spine. Once the complete paths missing part on the spine is calculated (l.42), they can be combined to build the path leading from found form to starting form (l.43). The rest is nothing more than converting this path which actually is a sequnce of generators on $\psl$ to a matrix (l.44), and then to decide which of the two orthogonal basis solves the form when multiplied with the path in matrix form (l.45-49).

For the algorithm listing of functions referenced in Algorithm.1 see Appendix, for the listing of \c{c}ark reduction algorithm see \cite{reduction}.

\section{Conclusions}

In this paper we present the software package InfoMod for representing and operating on binary quadratic forms. It consists of a library identically implemented in Matlab and Python\footnote{source code available at \url{https://github.com/hayral/infomod}}, which can compute reduction of binary quadratic forms and solve the representation problem. It also contains an interactive visualization application available for web \footnote{ \url{http://math.gsu.edu.tr/azeytin/infomod/node/3}} and mobile \footnote{ \url{https://play.google.com/store/apps/details?id=air.com.hkn.infomod}} platforms, which allows to observe the intricate patterns of forms' properties plotted according to their natural topology, and practically compute the reductions/representation problem solutions without writing code to invoke the libraries.
 
\bibliographystyle{plain}
\bibliography{infomod-bib}

\begin{thebibliography}{10}

\bibitem{Buchmann}
Ingrid Biehl and Johannes Buchmann.
\newblock An analysis of the reduction algorithms for binary quadratic forms.
\newblock In {\em Voronoi's Impact on Modern Science}, pages 71--98, 1997.

\bibitem{bqf/vollmer}
Johannes Buchmann and Ulrich Vollmer.
\newblock {\em Binary quadratic forms:An algorithmic approach}, volume~20 of
  {\em Algorithms and Computation in Mathematics}.
\newblock Springer, Berlin, 2007.

\bibitem{buell/bqf}
Duncan~A. {Buell}.
\newblock {\em {Binary quadratic forms. Classical theory and modern
  computations.}}
\newblock New York, NY etc.: Springer-Verlag, 1989.

\bibitem{cohen/diaz/olivier}
Henri {Cohen}, Francisco {Diaz y Diaz}, and Michel {Olivier}.
\newblock {Subexponential algorithms for class group and unit computations.}
\newblock {\em {J. Symb. Comput.}}, 24(3-4):433--441, 1997.

\bibitem{cox/primes/of/the/form}
David~A. Cox.
\newblock {\em Primes of the form {$x^2 + ny^2$}}.
\newblock Pure and Applied Mathematics (Hoboken). John Wiley \& Sons, Inc.,
  Hoboken, NJ, second edition, 2013.
\newblock Fermat, class field theory, and complex multiplication.

\bibitem{crawford2006actionscript}
Stewart Crawford and Elizabeth Boese.
\newblock Actionscript: a gentle introduction to programming.
\newblock {\em Journal of Computing Sciences in Colleges}, 21(3):156--168,
  2006.

\bibitem{farb/margalit/MCG}
Benson {Farb} and Dan {Margalit}.
\newblock {\em {A primer on mapping class groups.}}
\newblock Princeton, NJ: Princeton University Press, 2011.

\bibitem{disquisitiones}
Carl~Friedrich Gauss.
\newblock {\em Disquisitiones arithmeticae}, volume 157.
\newblock Yale University Press, 1966.

\bibitem{jacobson/williams}
M.~J. Jacobson, Jr. and H.~C. Williams.
\newblock Modular arithmetic on elements of small norm in quadratic fields.
\newblock {\em Des. Codes Cryptogr.}, 27(1-2):93--110, 2002.
\newblock Special issue in honour of Ronald C. Mullin, Part II.

\bibitem{jacobson/williams/solving/pell/equation}
Michael~J. Jacobson, Jr. and Hugh~C. Williams.
\newblock {\em Solving the {P}ell equation}.
\newblock CMS Books in Mathematics/Ouvrages de Math\'ematiques de la SMC.
  Springer, New York, 2009.

\bibitem{jacobson/phd/thesis}
M.J. Jacobson~Jr.
\newblock {\em Subexponential Class Group Computation in Quadratic Orders}.
\newblock PhD thesis, Technische Universit\"{a}t Darmstadt, 1999.

\bibitem{lagarias/computational/complexity/of/pell}
J.C. {Lagarias}.
\newblock {On the computational complexity of determining the solvability or
  unsolvability of the equation $X^2-DY^2=-1$.}
\newblock {\em {Trans. Am. Math. Soc.}}, 260:485--508, 1980.

\bibitem{lagarias/worst/case/complexity/bound}
J.C. {Lagarias}.
\newblock {Worst-case complexity bounds for algorithms in the theory of
  integral quadratic forms.}
\newblock {\em {J. Algorithms}}, 1:142--186, 1980.

\bibitem{lagrange}
J.~L. {Lagrange}.
\newblock {\"Uber die L\"osung der unbestimmten Probleme zweiten Grades (1768).
  Aus dem Franz\"osischen \"ubersetzt und herausgegeben von {\it E. Netto}.}
\newblock {Leipzig: Wilhelm Engelmann. 131 S. $8^\circ$ (Ostwalds Klassiker Nr.
  146) (1904).}, 1904.

\bibitem{sawilla/silvester/williams}
R.E. {Sawilla}, A.K. {Silvester}, and H.C. {Williams}.
\newblock {A new look at an old equation.}
\newblock In {\em {Algorithmic number theory. 8th international symposium,
  ANTS-VIII Banff, Canada, May 17--22, 2008 Proceedings}}, pages 37--59.
  Berlin: Springer, 2008.

\bibitem{schonhage}
Arnold Sch{\"o}nhage.
\newblock Fast reduction and composition of binary quadratic forms.
\newblock In {\em Proceedings of the 1991 international symposium on Symbolic
  and algebraic computation}, pages 128--133. ACM, 1991.

\bibitem{shanks/class/number}
Daniel {Shanks}.
\newblock {Class number, a theory of factorization, and genera.}
\newblock {1969 Number Theory Institute, Proc. Sympos. Pure Math. 20, 415-440
  (1971).}, 1971.

\bibitem{UZD}
A.~M. Uluda\u{g}, A.~Zeytin, and M.~Durmu\c{s}.
\newblock Binary quadratic forms as dessins.
\newblock 2016.
\newblock to appear in J. Th\'{e}or. Nombres Bordeaux.

\bibitem{zagier/zetafunktionen/quadratische/zahlkorper}
Don~B. Zagier.
\newblock {\em {Zetafunktionen und quadratische K\"orper. Eine Einf\"uhrung in
  die h\"ohere Zahlentheorie}}.
\newblock {}, 1981.

\bibitem{reduction}
A.~Zeytin.
\newblock On reduction theory of binary quadratic forms.
\newblock {\em {Publ. Math. Debrecen}}, 89:203--221, 2016.

\end{thebibliography}

% Appendix
\appendix
\section*{APPENDIX}

In this section you can find the algorithm listings for the functions invoked by Algorithm 1: Representation Problem Solver.

\begin{algorithm}
\caption{Generate All Spinal Forms}

\begin{algorithmic}[1] 
\Require $f_{spinal}$ is a form on the spine
\Function{RevolveAroundSpine}{$f_{spinal}$} 
\State $e_1 \gets (1,0) , e_2 \gets (0,1)$ \Comment $e_1,e_2$ : standard oriented basis of $\ZZ^2$
\State $ L \gets \begin{pmatrix} 1&-1\\ 1&0 \end{pmatrix} , S \gets \begin{pmatrix} 0&-1\\ 1&0 \end{pmatrix}  $ \Comment{$\psl$ generators}
\State $list \gets \emptyset$ %\{f_{spinal}\}$
\State $f_{next} \gets f_{spinal}$
\Repeat
\State $f_{next} \gets f_{next} \cdot SL$ \Comment  $\psl$ group action
\If{$f_{next}(e_1) \times f_{next}(e_2) > 0$}
\State $f_{next} \gets f_{next} \cdot L$ \Comment{multiply again to obtain $L^2$}
\EndIf
\State $list \gets list \cup \{f_{next}\}$
\Until{$f_{next} = f_{spinal}$} \Comment{repeat until we are back to starting form}
\State \textbf{return} $list$
\EndFunction
\end{algorithmic}
\end{algorithm}

\begin{algorithm}
\caption{Path Between Two Forms on Spine as a Word}

\begin{algorithmic}[1] 
\Require $f_{from}$ and $f_{to}$ are forms on spine
\Function{PathOnSpine}{$f_{from},f_{to}$} 

\State $e_1 \gets (1,0) , e_2 \gets (0,1)$ \Comment $e_1,e_2$ : standard oriented basis of $\ZZ^2$
\State $ L \gets \begin{pmatrix} 1&-1\\ 1&0 \end{pmatrix} , S \gets \begin{pmatrix} 0&-1\\ 1&0 \end{pmatrix}  $ \Comment{$\psl$ generators}
\State $path \gets \varepsilon$ \Comment empty string
\State $f_{next} \gets f_{from}$
\Repeat

\State $f_{next} \gets f_{next} \cdot L$ \Comment  $\psl$ group action
\State $path \gets path + "L"$ \Comment  + : string concatenation
\If{$f_{next}(e_1) \times f_{next}(e_2) > 0$}
\State $f_{next} \gets f_{next} \cdot L$ \Comment{multiply again to obtain $L^2$}
\State $path \gets path + "L"$
\EndIf

\If{$f_{next} = f_{to}$}
	\State \textbf{break} \Comment leave the loop early
\EndIf

\State $f_{next} \gets f_{next} \cdot S$ 
\State $path \gets path + "S"$

\Until{$f_{next} = f_{to}$} \Comment{repeat until we arrive at $f_{to}$}
\State \textbf{return} $path$

\EndFunction
\end{algorithmic}
\end{algorithm}

\begin{algorithm}
\caption{Convert a sequence of $\psl$ generators represented as symbols from alphabet $\{S,L\}$ to a $2 \times 2$ matrix}

\begin{algorithmic}[1] 
\Require $path$ is a sequence of symbols from alphabet $\{S,L\}$
\Function{PathToMatrix}{$path$} 

\State $ L \gets \begin{pmatrix} 1&-1\\ 1&0 \end{pmatrix} , S \gets \begin{pmatrix} 0&-1\\ 1&0 \end{pmatrix} , M \gets \begin{pmatrix} 1&0\\ 0&1 \end{pmatrix}    $ \Comment{$L,S : \psl$ generators, $M :$ identity}
%\State $ M \gets \begin{pmatrix} 1&0\\ 0&1 \end{pmatrix}  $

\For{$i \gets 1..$ \Call{Len}{$path$}}
\If{$path[i] = "L"$}
\State $M \gets M \cdot L$
\Else
\State $M \gets M \cdot S$
\EndIf
\EndFor
\State \textbf{return} $M$

\EndFunction
\end{algorithmic}
\end{algorithm}

\end{document}